\documentclass[12pt,thmsa]{article}
\usepackage{amsmath}
\usepackage{graphicx}
\usepackage{amsfonts}
\usepackage{amssymb}
\usepackage{mathrsfs}

\newtheorem{thm}{Theorem}

\begin{document}

\begin{center}
{\Large \textbf{Asymptotic normality of an estimator of kernel-based conditional mean dependence measure}}

\bigskip

Terence Kevin  MANFOUMBI DJONGUET and  Guy Martial  NKIET 

\bigskip

\textsuperscript{}URMI, Universit\'{e} des Sciences et Techniques de Masuku,  Franceville, Gabon.

\bigskip

E-mail adresses : tkmpro95@gmail.com,    guymartial.nkiet@univ-masuku.org.

\bigskip
\end{center}

\noindent\textbf{Abstract.} We propose an estimator of   the kernel-based conditional mean dependence measure  obtained from an appropriate  modification of a naive   estimator based on usual empirical estimators. We then get asymptotic normality of this estimator both under conditional mean  independence hypothesis and under the alternative hypothesis. A  new  test for conditional mean independence   of  random variables valued into Hilbert  spaces is then  introduced.

\bigskip

\noindent\textbf{AMS 1991 subject classifications: }62E20, 46E22.

\noindent\textbf{Key words:} Asymptotic normality;  Kernel method;   Kernel-based conditional  dependence; Reproducing kernel Hilbert space;  Functional data analysis.
\section{Introduction}

\noindent Conditional mean dependence is a statistical property that is important to evaluate for given random variables. Indeed, many regression analysis problems consist in modeling conditional mean of a response variable $Y$ given a predictor variable $X$ using either linear models or nonparametric models. Such  modeling approaches  are  in fact not relevant in case of conditional mean independence of the involved variables. That is why testing whether the predictor  has  a contribution  to the mean of the response is of a great interest. However, there exists just a few works dealing with the problem of testing for conditional mean independence between random variables.  It was investigated in Shao and Zhang (2014)  by using the so-called martingale difference divergence (MDD) for $Y\in\mathbb{R}$ and $X\in\mathbb{R}^q$. Later, a generalization of MDD was introduced inPark et al.  (2015)  in order to deal with the case of  $Y\in\mathbb{R}^p$ and $X\in\mathbb{R}^q$, and Lee et al.  (2020)  proposed functional martingale difference divergence (FMDD) which extended MDD to the case where $X$ and $Y$ are functional variables. The case of high-dimensional setting was tackled in Zhang  et al.  (2018) .   Recently, Lai  et al.  (2021)   introduced the kernel-based  conditional mean dependence measure (KCMD) by means of which  a test for conditional mean independence was constructed.  This test is based on an unbiased estimator of KCMD which has the form of a U-statistic  with   asymptotic distribution under null hypothesis  equal to   an infinite sum of distributions. This last property is a drawback that forced Lai  et al.  (2021)  to resort to a wild boostrap method for performing the test.  Faced with a similar problem with a maximal mean discrepancy (MMD) estimator, Magikusa and Naito  (2020)  adopted an approach permitting to obtain asymptotic normality for a proposed estimator both under the null hypothesis and under the alternative. This approach was also  used later in Balogoun et al.   (2021)  for the case of generalized maximal mean discrepancy (GMMD). In this paper  we tackle this approach consisting in making an appropriate modification on a naive estimator of KCMD. We then obtain asymptotic normality for the resulting estimator under the conditional mean independence hypothesis. This allows to propose a test for conditional mean independence of random variables with values into Hilbert spaces and that can, therefore, be used on functional data. The rest of the paper is organized as follows.     The KCMD  is recalled in Section 2, and Section 3 is devoted to its estimation by a modification of the naive estimator,  and to the main results. All the proofs are postponed in Section 4.

\section{KCMD and   conditional mean independence}
Let $X$ and $Y$ be two  random variables  defined on a probability space $(\Omega,\mathscr{A},P)$ and   taking values in separable Hilbert spaces $\mathcal{X}$ and $\mathcal{Y}$ respectively; it is assumed that  $\mathbb{E} \left(\Vert Y\Vert_\mathcal{Y} ^2\right) < +\infty$, where $\Vert\cdot\Vert_\mathcal{Z}$ denotes the norm associated with the inner product  $\langle \cdots,\cdot  \rangle_\mathcal{Z}$ of the Hilbert space $\mathcal{Z}$ . In order to test for conditional mean independence, that is testing for the hypothesis 
$$ \mathscr{H}_0\, : \, \mathbb{E} \left(Y|X\right)=\mathbb{E} \left(Y\right) \mbox{  almost surely     } $$
versus
$$\mathscr{H}_1\, :\, P\big( \mathbb{E} \left(Y|X\right)=\mathbb{E} \left(Y\right) \big) < 1,$$
where $\mathbb{E} \left(Y|X\right)$ denotes conditional expectation,  Lai  et al.  (2021)  introduced the Kernel Conditional Mean Independence measure  (KCMD). Let us consider a reproducing kernel Hilbert space    $\mathcal{H}$  of functions  from  $\mathcal{X}$ to $\mathbb{R}$  with associated   kernel    $K\,:\,\mathcal{X}^2\rightarrow\mathbb{R}$ which is  a symmetric function such that, for any $f\in \mathcal{H}$ and any $x\in\mathcal{X}$, one has $K(x,\cdot)\in\mathcal{H}$ and $f(x)=\langle K(x,\cdot),f\rangle_{\mathcal{H}}$   (see  Berlinet and Thaomas-Agnan (2004)).   Throughout this paper, we assume that $K$    satisfies  the following condition:

\bigskip

\noindent $(\mathscr{C}_1):$  $\Vert  K\Vert_{\infty}:= \sup\limits_{(x,y)\in \mathcal{X}^2} K(x,y)<+ \infty $;
 
\bigskip
\noindent then  the kernel mean embedding  $m_X:=\mathbb{E}\left(K(X,\cdot)\right)$ exists. KCMD is the measure given by 
\begin{equation}\label{KCMD}
\begin{aligned}
\mbox{KCMD(Y,X)}  & = \bigg\Vert \mathbb{E} \bigg(Y \otimes K(X,\cdot)\bigg) -    \mu \otimes m_X \bigg\Vert_{\textrm{HS}} ^2,
\end{aligned}
\end{equation}
where $\mu=\mathbb{E}(Y)$, the tensor product $\otimes$ is such that, for any $(y,f)\in \mathcal{Y}\times\mathcal{H}$, $y\otimes f$ is the linear operator   defined by  $(y\otimes f)(t)=\langle y,t\rangle_{\mathcal{Y} }\,f$  for any $t\in\mathcal{Y} $, and $ \Vert \cdot\Vert_{\textrm{HS}}$ denotes the Hilbert-Schmidt norm of operators. As demonstrated in Lai  et al.  (2021) , when the kernel $K$ is   characteristic, then the null hypothesis $\mathscr{H}_0$ holds if, and only if,    $\mbox{KCMD(Y,X)} =0$. So, a test for conditional mean independence can be achieved by using an estimator of  $\mbox{KCMD(Y,X)} =0$ as test statistic. An unbiased estimator, based on a i.i.d. sample   $\{(X_i,Y_i)\}_{1 \leq i \leq n}$   of $(X,Y)$,   was defined in Lai  et al.  (2021)  as:
\[
\mbox{KCMD}_n(Y,X) = \frac{1}{n(n-3)} \sum_{i \neq j} C_{ij}D_{ij},
\]
where
\[
c_{ij} = \left\{
\begin{array}{lcl}
K(X_i,Y_j)& &\textrm{if }  i \neq j  \\
 0 & &\textrm{if }  i=j 
\end{array} 
\right.,    \qquad     
d_{ij} = \left\{
\begin{array}{lcl} 
\langle Y_i,Y_j\rangle_\mathcal{Y}& &\textrm{if }  i \neq j \\ 
 0 & &\textrm{if }  i=j 
\end{array}  
\right.,
\]
\[
c_{i\cdot} = \frac{1}{n-2} \sum_{j=1}^{n}c_{ij},   \qquad     c_{\cdot j} = \frac{1}{n-2} \sum_{i=1}^{n}c_{ij}, \qquad     c_{\cdot\cdot} = \frac{1}{(n-1)(n-2)} \sum_{i,j=1}^{n}c_{ij},
\]
\[
d_{i\cdot} = \frac{1}{n-2} \sum_{j=1}^{n}d_{ij},   \qquad     d_{\cdot j} = \frac{1}{n-2} \sum_{i=1}^{n}d_{ij}, \qquad     d_{\cdot\cdot} = \frac{1}{(n-1)(n-2)} \sum_{i,j=1}^{n}d_{ij},
\]
\[
C_{ij} = \left\{
\begin{array}{lcl} 
c_{ij}-c_{i\cdot}-c_{\cdot j}+c_{\cdot\cdot}  & &\textrm{if }  i \neq j  \\
 0 & &\textrm{if }  i=j 
\end{array} \right.,
\]
\[     
D_{ij} =  \left\{
\begin{array}{lcl} 
d_{ij}-d_{i\cdot}-d_{\cdot j}+d_{\cdot\cdot}  & &\textrm{if }  i \neq j  \\
 0 & &\textrm{if }  i=j 
\end{array} \right..
\]
They derived the asymptotic distribution under null hypothesis of this statistic and proved that, under $\mathscr{H}_0$,  $n\mbox{KCMD}_n(Y,X) $ converges in distribution, as $n\rightarrow +\infty$, to $\sum_{m=1}^{+\infty}\gamma_m\left(N_m^2-1\right)$, where  $N_m$s are i.i.d. standard normal distributed random variables and $(\gamma)_{m\geq 1}$ is a sequence of eigenvalues of a suitable positive autoadjoint operator.  This limiting distribution can not be used to compute critical values for performing the test since the $\gamma_m$s are unkown,  and since it is an infinite sum of distributions. That is why Lai  et al.  (2021)  proposed a wild  bootstrap procedure to approxmate  the asymptotic null distribution.  As one knows, bootstrap procedures have the disadvantage of leading to rather high computation times, that is why it is preferable to obtain asymptotic normality of the test statistic. Following an approach introduced inMagikusa and Naito  (2020)   and also tackled in Balogoun et al.   (2021) , we propose in this paper to modify a naive estimator of KCMD$(Y,X)$ in order to get asymptotic normality under $\mathscr{H}_0$ and to use this result for performing the test.

\section{Modification of KCMD and asymptotic normality}

Replacing each expectation in \eqref{KCMD} by its empirical counterpart leads  to the simple estimator of KCMD  given by
\[
\widehat{\mbox{KCMD}}_n  =\bigg \Vert \frac{1}{n} \sum_{i=1}^{n} Y_i \otimes K(X_i,\cdot) - \bigg(\frac{1}{n} \sum_{i=1}^{n} Y_i \bigg) \otimes \bigg(\frac{1}{n} \sum_{i=1}^{n} K(X_i,\cdot)\bigg) \bigg\Vert_{\textrm{HS}} ^2
\]
and which can be expanded as
\begin{eqnarray}\label{esti}
\widehat{\mbox{KCMD}}_n  &=&\frac{1}{n^2}\sum_{i,j=1}^{n}  \langle Y_i , Y_j \rangle_\mathcal{Y} K(X_i,X_j) + \frac{1}{n^4}\sum_{i,j,q,r=1}^{n}  \langle Y_q , Y_r \rangle_\mathcal{Y} K(X_i,X_j) \nonumber\\
& &- \frac{2}{n^3} \sum_{i,j,q=1}^{n}  \langle Y_i , Y_q \rangle_\mathcal{Y} K(X_i,X_j)
\end{eqnarray}
by using properties of $\otimes$ and reproducing property of $K$.  We propose another estimator of KCMD$(Y,X)$ obtained from a modification of $\widehat{\mbox{KCMD}}_n$. This modification just consists to introduce a weight in the croos-product term of \eqref{esti}. Let  $\left\{w_{i,n}(\gamma)\right\}_{1\leq i\leq n}$ be a triangular array  of positive real numbers depending on a parameter $\gamma\in]0,1[$. We consider the  estimator  $\widehat{K}_{n,\gamma} $ of  KCMD$(Y,X)$ given by:
\begin{eqnarray*}\label{esti2}
\widehat{K}_{n,\gamma}   &=&\frac{1}{n^2}\sum_{i,j=1}^{n}  \langle Y_i , Y_j \rangle_\mathcal{Y} K(X_i,X_j) + \frac{1}{n^4}\sum_{i,j,q,r=1}^{n}  \langle Y_q , Y_r \rangle_\mathcal{Y} K(X_i,X_j) \nonumber\\
& &- \frac{2}{n^3} \sum_{i,j,q=1}^{n}w_{i,n}(\gamma)  \langle Y_i , Y_q \rangle_\mathcal{Y} K(X_i,X_j),
\end{eqnarray*}
and we take it as test statistic. For obtaining its asymptotic normality, we suppose that   the sequence of weights that is used satisfy the following conditions:

\bigskip

\noindent $(\mathscr{C}_2):$ There exists  a strictly positive real number $ \tau $  and an  integer $ n_0 $ such that for all $n> n_0$:
$$n\left\vert\frac{1}{n}\sum_{i=1}^{n}w_{i,n}(\gamma)-1\right\vert\leq \tau.$$

\noindent $(\mathscr{C}_3):$ There exists  $C>0$ such that $\max\limits_{1\leq i \leq n}w_{i,n}(\gamma)<C$ for all  $n\in\mathbb{N}^\ast$ and $ \gamma \in]0, 1]$.

\noindent$(\mathscr{C}_4):$ For any   $ \gamma \in]0, 1[$,
$\lim\limits_{n\rightarrow +\infty}\frac{1}{n}\sum_{i=1}^{n}w^2_{i,n}(\gamma)=w^2(\gamma)> 1$.

 \bigskip

\noindent Such sequence was first introduced in \cite{Ahmad1993} where an example defined as   $w_{i,n}(\gamma)=1+(-1)^i\,\gamma$ was given.  For this example, one has  $C=2$, $w^2(\gamma)=1+\gamma^2$ and $\tau$ is any positive real number. Another example is   $w_{i,n}(\gamma)=1+\sin(i\pi\gamma)$ which corresponds to $\tau=1/\vert\sin(\pi\gamma/2)\vert$, $C=2$ and $w^2(\gamma)=3/2$.
\noindent Putting $
\eta=\mathbb{E}\left( Y  \otimes K(X,\cdot)\right)$ and $ \nu=\mu\otimes m_X$, 
and considering the functions  $\mathcal{U}$ and $\mathcal{V}$  from $\mathcal{X}\times\mathcal{Y}$ to $\mathbb{R}$ defined as
\begin{equation*}\label{fonc1}
\begin{aligned}
 \mathcal{U} (x,y) &= \big\langle  y  \otimes K(x,\cdot) - \eta ,  \eta \big\rangle_{\textrm{HS}}+\big\langle y \otimes m_X + \mu \otimes  K(x,\cdot) - 2 \nu ,  \nu - \eta \big\rangle_{\textrm{HS}},\\
\mathcal{V}(x,y) &= \big\langle  y  \otimes K(x,\cdot)- \eta ,  \nu \big\rangle_{\textrm{HS}},
\end{aligned}
\end{equation*}
where $\big\langle   \cdot , \cdot  \big\rangle_{\textrm{HS}}$  denotes the Hilbert-Smidt inner product,  we have:

\bigskip

\begin{thm}{}\label{thm}
Assume  conditions  $(\mathscr{C}_1)$ to  $(\mathscr{C}_4)$. Then   $$ \sqrt{n}\Big (\widehat{K}_{n,\gamma} - \textrm{ KCMD}(Y,X) \Big ) \stackrel{\mathscr D}{\longrightarrow} \mathcal N (0, \sigma ^2_\gamma),$$ as $n\rightarrow +\infty$,  where 
\begin{equation*}
\begin{aligned}
\sigma ^2_\gamma   & =  4 Var\left(\mathcal{U}(X_1,Y_1)\right)+4 w^2(\gamma)Var\left(\mathcal{V}(X_1,Y_1)\right) -8Cov\left(\mathcal{U}(X_1,Y_1),\mathcal{V}(X_1,Y_1)\right).
\end{aligned}
\end{equation*}
\end{thm}

\bigskip

\noindent This theorem gives asymptotic normality both under $\mathscr{H}_0$ and under $\mathscr{H}_1$. Under $\mathscr{H}_0$, we have  KCMD$(Y,X)=0$, which is equivalent to $\eta=\nu$ and implies that $ \mathcal{U} = \mathcal{V} $;  then, $ \sqrt{n} \widehat{K}_{n,\gamma}  \stackrel{\mathscr D}{\longrightarrow} \mathcal N (0, \sigma ^2_\gamma)$, as $n\rightarrow +\infty$, with
$\sigma ^2_\gamma   =  4 \left(w(\gamma) ^2- 1\right)Var\left( \big\langle  Y_1 \otimes K(X_1,\cdot),  \eta \big\rangle_{\textrm{HS}}\right)$.  This variance is unknown since  its depends on $\nu$. So, for performing the test  we have to estimate it. We consider the estimator  
\[
\widehat{\sigma}_\gamma^2=4 (w(\gamma) ^2- 1)\widehat{\alpha}^2 ,
\] 
where
\begin{equation*}\label{sigma}
\begin{aligned}
\widehat{\alpha}^2 & =  \frac{1}{n}\sum_{i=1}^{n}  \Big ( \frac{1}{n}\sum_{j=1}^{n}  \langle Y_i , Y_j \rangle_\mathcal{Y} K(X_i,X_j)  - \frac{1}{n^2}\sum_{m=1}^{n} \sum_{p=1}^{n}  \langle Y_m , Y_p \rangle_\mathcal{Y} K(X_m,X_p)  \Big) ^2,
\end{aligned}
\end{equation*}
and we have:

\bigskip

\begin{thm}{}\label{thm2}
Assume  conditions  $(\mathscr{C}_1)$ to  $(\mathscr{C}_4)$. Then, under $\mathscr{H}_0$,
\[
\sqrt{n}\frac{\widehat{K}_{n,\gamma}}{\widehat{\sigma}_\gamma}  \stackrel{\mathscr D}{\longrightarrow} \mathcal N (0, 1),
\] 
as $n\rightarrow +\infty$.
\end{thm}

\bigskip

\noindent This theorem allows to achieve the test in practice. The null hypothesis $\mathscr{H}_0$ is to be rejected when $\sqrt{n} \widehat{K}_{n,\gamma}> \widehat{\sigma}_\gamma  \Phi^{-1}(1-\alpha)$, where $\alpha$ is the chosen significance level and $\Phi$ is the cumulative distribution function of the standard normal distribution.

\bigskip

\noindent\textbf{Remark 1.}
This test can be applied on functional data corresponding, for instance,  to the case where the $X_i$s and the $Y_i$s are random functions belonging in $L^2([0,1])$ and observed on points $t_1,\cdots,t_r$ and $s_1,\cdots,s_q$, respectively, of   fine grids  in $[0,1]$ such that $t_1=s_1=0$ and $t_r=s_q=1$. In this case, one has
\[
\langle Y_i , Y_j \rangle_\mathcal{Y}=\int_0^1Y_i(s)\,Y_j(s)\,ds,
\]
what can be approximated  by using trapezoidal rule so as to obtain
\begin{equation}\label{approx1}
\langle Y_i , Y_j \rangle_\mathcal{Y}\simeq\sum_{m=1}^{q-1}\frac{s_{m+1}-s_m}{2}\left(Y_i(s_m)\,Y_j(s_m)+Y_i(s_{m+1})\,Y_j(s_{m+1})\right).
\end{equation}
If the gaussian kernel is used, one has
\[
K(X_i ,X_j )=\exp\left(-\omega^2\Vert X_i-X_j\Vert^2_\mathcal{X}\right)=\exp\left(-\omega^2\int_0^1\left(X_i(t)-X_j(t)\right)^2\,dt\right),
\]
where $\omega>0$, and this term  can also  be approximated  by using trapezoidal rule:
\begin{equation}\label{approx2}
K(X_i ,X_j )\simeq\exp\left(-\omega^2\sum_{m=1}^{r-1}\frac{t_{m+1}-t_m}{2}\left(\left(X_i(t_m)-X_j(t_m)\right)^2+\left(X_i(t_{m+1})-X_j(t_{m+1})\right)^2\right)\right).
\end{equation}
Then,   $\widehat{K}_{n,\gamma}$ and $\widehat{\alpha}^2$  are to be computed by using \eqref{approx1} and \eqref{approx2}.

\section{Proofs}

\subsection{Proof of Theorem \ref{thm}}

Putting  $ \overline{Y}_n = \frac{1}{n} \sum_{i=1}^{n} Y_i$, $\overline{K}_n = \frac{1}{n} \sum_{i=1}^{n} K(X_i,\cdot)$ and $\overline{Y \otimes K}^n = \frac{1}{n} \sum_{i=1}^{n}Y_i\otimes  K(X_i,\cdot)$, we have
\begin{equation*}
\begin{aligned}
&\sqrt{n}\Big (\widehat{K}_{n,\gamma}  -  \textrm{KCMD}(Y,X)  \Big ) \\
& = \sqrt{n}\Bigg [  \bigg\Vert \overline{Y \otimes K}^n - \eta \bigg\Vert_{\textrm{HS}} ^2 + 2 \; \big\langle \overline{Y \otimes K}^n ,  \eta  \big\rangle_{\textrm{HS}} - \Vert \eta \Vert_{\textrm{HS}} ^2  + \bigg\Vert \overline{Y}_n \otimes \overline{K}_n - \nu \bigg\Vert_{\textrm{HS}} ^2 + 2 \; \big\langle \overline{Y}_n \otimes \overline{K}_n ,  \nu  \big\rangle_{\textrm{HS}} \\
& \qquad  - \Vert \nu \Vert_{\textrm{HS}} ^2 - \frac{2}{n} \sum_{i=1}^{n} \left(w_{i,n}(\gamma) - 1 \right) \;  \big\langle  Y_i \otimes K(X_i,\cdot) , \overline{Y}_n \otimes \overline{K}_n  - \nu \big\rangle_{\textrm{HS}}\\
& \qquad- \frac{2}{n} \sum_{i=1}^{n}w_{i,n}(\gamma) \;  \big\langle  Y_i \otimes K(X_i,\cdot) ,  \nu \big\rangle_{\textrm{HS}}\\
& \qquad  - \Vert \eta - \nu \Vert_{\textrm{HS}} ^2 - 2  \big\langle \overline{Y \otimes K}^n - \eta ; \overline{Y}_n \otimes \overline{K}_n  - \nu \big\rangle_{\textrm{HS}} - 2 \; \big\langle \overline{Y}_n \otimes \overline{K}_n ,  \eta  \big\rangle_{\textrm{HS}}\\
& \qquad  - \frac{2}{n} \sum_{i=1}^{n} \left(w_{i,n}(\gamma) - 1 \right) \;  \big\langle \eta ,  \nu \big\rangle_{\textrm{HS}} +\frac{2}{n} \sum_{i=1}^{n}w_{i,n}(\gamma) \;  \big\langle \eta ,  \nu \big\rangle_{\textrm{HS}} \Bigg ] \\
& = A_n + B_n + C_n + D_n,
\end{aligned}
\end{equation*}
where
\begin{equation*}\label{test}
\begin{aligned}
A_n &=n^{-1/2} \bigg ( \bigg\Vert \sqrt{n}\big(\overline{Y \otimes K}^n - \eta \big)\bigg\Vert_{\textrm{HS}} ^2+ \bigg\Vert \sqrt{n}\big(\overline{Y}_n \otimes \overline{K}_n - \nu \big)\bigg\Vert_{\textrm{HS}} ^2\bigg)\\
& = n^{-1/2} \bigg( \bigg\Vert \sqrt{n}\big(\overline{Y \otimes K}^n - \eta \big)\Vert_{\textrm{HS}} ^2+ \bigg\Vert \left(\sqrt{n}(\overline{Y}_n - \mu)\right)\otimes\left(\overline{K}_n - m_X\right) \\
& \qquad \qquad +\left(\sqrt{n}(\overline{Y}_n - \mu)\right)\otimes m_X + \mu \otimes\left(\sqrt{n}(\overline{K}_n - m_X)\right) \bigg\Vert_{\textrm{HS}} ^2\bigg)\\
\end{aligned}
\end{equation*}
\begin{equation*}
\begin{aligned}
B_n 
& = -2 \bigg(  \;  \Big\langle \frac{1}{n} \sum_{i=1}^{n} \left(w_{i,n}(\gamma) - 1 \right) Y_i \otimes K(X_i,\cdot) , \sqrt{n} \big(\overline{Y}_n \otimes \overline{K}_n  - \nu \big) \Big\rangle_{\textrm{HS}} \\
&\qquad+ \frac{\sqrt{n}}{n} \sum_{i=1}^{n} \left(w_{i,n}(\gamma) - 1 \right) \;  \big\langle \eta ,  \nu \big\rangle_{\textrm{HS}}\bigg),\\
\end{aligned}
\end{equation*}
\begin{equation*}
\begin{aligned}
C_n &= - 2  \Big\langle \sqrt{n} \big(\overline{Y \otimes K}^n - \eta \big), \overline{Y}_n \otimes \overline{K}_n  - \nu \Big\rangle_{\textrm{HS}}+ 2\big\langle \sqrt{n}(\overline{Y}_n - \mu) \otimes (\overline{K}_n - m_X),  \nu - \eta \big\rangle_{\textrm{HS}},  \\
\end{aligned}
\end{equation*}
\begin{equation*}
\begin{aligned}
D_n &= \frac{2}{ \sqrt{n}}\sum_{i=1}^{n} \mathcal{U} (X_i,Y_i) - w_{i,n}(\gamma)\,\mathcal{V} (X_i,Y_i). \nonumber\\
\end{aligned}
\end{equation*}
The central limit theorem ensures that $\sqrt{n}\big(\overline{Y \otimes K}^n - \eta \big)$, $\sqrt{n}\big(\overline{Y}_n - \mu \big)$ and $\sqrt{n}\big(\overline{K}_n - m_X \big)$ converge in distribution to random variables having normal distributions as $n\rightarrow +\infty$. Moreover, by the law of large numbers $\overline{K}_n - m_X$ converges in probability to 0 as $n\rightarrow +\infty$.  Then, by the continuous mapping theorem we deduce that $A_n=o_p(1)$. Concerning $B_n$, we get by the Cauchy-Schwarz inequality 
\begin{equation*}
\begin{aligned}
|B_n| & \leq 2 \Bigg [\bigg\Vert \frac{1}{n} \sum_{i=1}^{n} \left(w_{i,n}(\gamma) - 1 \right) Y_i \otimes K(X_i,\cdot)\bigg\Vert_{\textrm{HS}}\,\bigg\Vert  \sqrt{n} \big(\overline{Y}_n \otimes \overline{K}_n  - \nu \big)  \bigg\Vert _{HS} \\
& \qquad \qquad + \sqrt{n} \Big |  \frac{1}{n} \sum_{i=1}^{n} \left(w_{i,n}(\gamma) - 1 \right)  \Big | \;\left\Vert \eta\right\Vert _{\textrm{HS}}\left\Vert \nu\right\Vert _{\textrm{HS}}\Bigg ]; \\
\end{aligned}
\end{equation*}
and since  $\sqrt{n} \big(\overline{Y}_n \otimes \overline{K}_n  - \nu \big)  = \sqrt{n} \big(\overline{Y}_n - \mu\big)\otimes \overline{K}_n  + \mu \otimes \sqrt{n} \big(\overline{K}_n - m_ X\big)$ and, under $(\mathscr{C}_2)$ for $n$ large enough,   $n\Big|\frac{1}{n}\sum_{i=1}^{n}w_{i,n}(\gamma) - 1 \Big| \leq\tau$,  it follows
\begin{equation*}
\begin{aligned}
|B_n| & \leq 2 \Bigg [\bigg\Vert \frac{1}{n} \sum_{i=1}^{n} \left(w_{i,n}(\gamma) - 1 \right) Y_i \otimes K(X_i,\cdot)\bigg\Vert_{\textrm{HS}} \Bigg (\bigg\Vert \sqrt{n} \big(\overline{Y}_n - \mu\big)\otimes \overline{K}_n \bigg\Vert _{HS} \\
& \qquad \qquad +\bigg\Vert    \mu \otimes \sqrt{n} \big(\overline{K}_n - m_ X\big)\bigg\Vert_{\textrm{HS}}  \Bigg)   +  \frac{ \tau}{\sqrt{n}}\left\Vert \eta\right\Vert _{\textrm{HS}}\left\Vert \nu\right\Vert _{\textrm{HS}}\Bigg ] \\
& = 2 \Bigg [\bigg\Vert \frac{1}{n} \sum_{i=1}^{n} \left(w_{i,n}(\gamma) - 1 \right) Y_i \otimes K(X_i,\cdot)\bigg\Vert_{\textrm{HS}} \Bigg (\bigg\Vert  \sqrt{n} \big(\overline{Y}_n - \mu\big)\bigg\Vert_{\mathcal{Y}}\;\bigg\Vert \overline{K}_n \bigg\Vert _{\mathcal{H}} \\
& \qquad \qquad +\left\Vert    \mu\right\Vert _{\mathcal{Y}} \;\bigg\Vert \sqrt{n} \big(\overline{K}_n - m_ X\big) \bigg\Vert _{\mathcal{H}} \Bigg)   +  \frac{ \tau}{\sqrt{n}}\left\Vert \eta\right\Vert _{\textrm{HS}}\left\Vert \nu\right\Vert _{\textrm{HS}}\Bigg ]. \\
\end{aligned}
\end{equation*}
By the reproducing property we obtain $$\left\Vert \overline{K}_n  \right\Vert_{\mathcal{H}} \leq \frac{1}{n}\sum_{i=1}^{n} \Vert K(X_i,\cdot) \Vert_{\mathcal{H}_X} \leq\left\Vert  K  \right\Vert_\infty^{1/2};$$  hence
\begin{equation*}
\begin{aligned}
|B_n| & \leq 2 \Bigg [\bigg\Vert \frac{1}{n} \sum_{i=1}^{n} \left(w_{i,n}(\gamma) - 1 \right) Y_i \otimes K(X_i,\cdot)\bigg\Vert_{\textrm{HS}} \Bigg (\bigg\Vert  \sqrt{n} \big(\overline{Y}_n - \mu\big)\bigg\Vert_{\mathcal{Y}}\;\left\Vert  K  \right\Vert_\infty^{1/2} \\
& \qquad \qquad +\left\Vert    \mu\right\Vert _{\mathcal{Y}} \;\bigg\Vert \sqrt{n} \big(\overline{K}_n - m_ X\big) \bigg\Vert _{\mathcal{H}} \Bigg)   +  \frac{ \tau}{\sqrt{n}}\left\Vert \eta\right\Vert _{\textrm{HS}}\left\Vert \nu\right\Vert _{\textrm{HS}}\Bigg ].
\end{aligned}
\end{equation*}
From Lemma 1 in Manfoumbi Djonguet et al. (2022)  we have $$\bigg\Vert\frac{1}{n} \sum_{i=1}^{n} \left(w_{i,n}(\gamma) - 1 \right) Y_i \otimes K(X_i,\cdot)\bigg\Vert_{\textrm{HS}} = o_p(1),$$ and from   the central limit theorem $\sqrt{n} \big(\overline{K}_n - m_ X\big)$ and $\sqrt{n} \big(\overline{Y}_n - \mu\big)$ converge in distribution as $n\rightarrow +\infty$. We then deduce from the preceding inequality that $B_n = o_p(1)$. Another use of   the Cauchy-Schwartz inequality yields:
\begin{equation*}
\begin{aligned}
|C_n| & \leq 2 \bigg\Vert  \sqrt{n} \big(\overline{Y \otimes k}^n - \eta \big)  \bigg\Vert _{\textrm{HS}}\;  \bigg\Vert  \overline{Y}_n \otimes \overline{K}_n  - \nu   \bigg\Vert _{\textrm{HS}}  \\
&\qquad+ 2\bigg\Vert  \sqrt{n}(\overline{Y}_n - \mu) \otimes (\overline{K}_n - m_X) \bigg\Vert _{\textrm{HS}}\; \bigg\Vert   \nu - \eta \bigg\Vert _{\textrm{HS}}\\
& = 2 \bigg\Vert  \sqrt{n} \big(\overline{Y \otimes k}^n - \eta \big)  \bigg\Vert _{\textrm{HS}}\;  \bigg\Vert  \overline{Y}_n \otimes \overline{K}_n  - \nu   \bigg\Vert _{\textrm{HS}}  \\
&\qquad+ 2\bigg\Vert  \sqrt{n}(\overline{Y}_n - \mu)\bigg\Vert _{\mathcal{Y}} \bigg\Vert \overline{K}_n - m_X \bigg\Vert _{\mathcal{H}}\; \bigg\Vert   \nu - \eta \bigg\Vert _{\textrm{HS}}.\\
\end{aligned}
\end{equation*}
As $n \rightarrow +\infty$, $\sqrt{n} \big(\overline{Y \otimes k}^n - \eta \big)$ and $ \sqrt{n}(\overline{Y}_n - \mu)$ converge in distribution to normal random variables, $\overline{K}_n$ and $\overline{Y}_n$ converge in probability to $m_X$ and $\mu$ respectively. Thus, by the continuous mapping theorem, $\overline{Y}_n \otimes \overline{K}_n$ converge in probability to $\nu$ as $n \rightarrow +\infty$, and the preceding inequality implies that $C_n=o_p(1)$. Finally, we got 
\begin{equation*}
\sqrt{n}\Big (\widehat{K}_{n,\gamma}  -  \textrm{KCMD}(Y,X)  \Big )= \frac{2}{ \sqrt{n}}\sum_{i=1}^{n} \mathcal{U} (X_i,Y_i) - w_{i,n}(\gamma)\,\mathcal{V} (X_i,Y_i) + o_p (1)  .
\end{equation*}
From Slutsky's  theorem,  $\sqrt{n}\Big (\widehat{K}_{n,\gamma}  -  \textrm{KCMD}(Y,X)  \Big )$ has the same limiting distribution than $D_n$.  Let us set  
\[
s_{n,\gamma}^2=\sum_{i=1}^{n}Var\Big(\mathcal{U}(X_i,Y_i)-w_{i,n}(\gamma)\mathcal{V}(X_i,Y_i)\Big).
\] 
 By similar arguments as in the proof of Theorem 1 inMagikusa and Naito  (2020)  we obtain that, for any $\varepsilon>0$,
\[
s_{n,\gamma}^{-2}\sum_{i=1}^{n}\int_{\{(x,y):| \mathcal{U}(x,y)-w_{i,n}(\gamma)\mathcal{V}(x,y) |>\varepsilon s_{n,\gamma}\}}^{} \bigg(\mathcal{U}(x,y)-w_{i,n}(\gamma)\mathcal{V}(x,y)\bigg)^2\,d\mathbb{P}_{XY}(x,y) 
\]
converges to $0$ as $n\rightarrow +\infty$. Therefore, by Section 1.9.3 in \cite{Serfling1980} we obtain that  
$\sqrt{n}s_{n,\gamma}^{-1}\frac{D_n}{2}	\stackrel{\mathscr{D}}{\rightarrow} \mathcal{N}\left(0,1\right)$. However,
\begin{eqnarray*}
\left(\frac{s_{n,\gamma}}{\sqrt{n}}\right)^2
&=&Var\left(\mathcal{U}(X_1,Y_1)\right)+\left(\frac{1}{n}\sum_{i=1}^{n}w^2_{i,n}(\gamma)\right)Var\left(\mathcal{V}(X_1,Y_1)\right)\nonumber\\
&&-2\left(\frac{1}{n}\sum_{i=1}^{n}w_{i,n}(\gamma)\right)Cov\left(\mathcal{U}(X_1,Y_1),\mathcal{V}(X_1,Y_1)\right),
\end{eqnarray*}
then, using $(\mathscr{C}_2)$ and  $(\mathscr{C}_4)$, we get
\[
\lim\limits_{n\rightarrow +\infty}\left(n^{-1}s_{n,\gamma}^2\right)=Var\left(\mathcal{U}(X_1,Y_1)\right)+ w^2(\gamma)Var\left(\mathcal{V}(X_1,Y_1)\right) -2Cov\left(\mathcal{U}(X_1,Y_1),\mathcal{V}(X_1,Y_1)\right)
\]
and, therefore, $D_n	\stackrel{\mathscr{D}}{\rightarrow} \mathcal{N}\left(0,\sigma_\gamma^2\right)$. 

\subsection{Proof of  Theorem \ref{thm2}}

It suffices to prove that $\widehat{\sigma}^2_\gamma$  converges in probability to $\sigma^2_\gamma$ as $n\rightarrow +\infty$, what  is obtained from the convergence in probability of $\widehat{\alpha}^2$  to  $Var\left( \big\langle  Y_1 \otimes K(X_1,\cdot),  \eta \big\rangle_{\textrm{HS}}\right)$.  From the definition of the Hilbert-Schmidt inner  product and the reproducing property of $K$ one can easily see that 
\[
\big\langle Y_i \otimes K(X_i,\cdot) , \overline{Y \otimes K}^n \big\rangle_{\textrm{HS}}=n^{-1} \sum_{i=1}^{n}\big\langle Y_i  , Y_j \big\rangle_{\mathcal{Y}}K(X_i,X_j)
\]
and, therefore, that
\begin{equation}\label{alpha}
\widehat{\alpha}^2 =\frac{1}{n}\sum_{i=1}^{n}   \big\langle  Y_i \otimes K(X_i,\cdot) , \overline{Y \otimes K}^n \big\rangle_{\textrm{HS}}^2  - \bigg(  \frac{1}{n}\sum_{i=1}^{n}   \big\langle  Y_i \otimes K(X_i,\cdot) , \overline{Y \otimes K}^n \big\rangle_{\textrm{HS}}   \bigg)^2.
\end{equation}
Noticing  that
\begin{equation}\label{diff}
\begin{aligned}
& \frac{1}{n}\sum_{i=1}^{n}   \big\langle  Y_i \otimes K(X_i,\cdot) , \overline{Y \otimes K}^n \big\rangle_{\textrm{HS}}^2 - \frac{1}{n}\sum_{i=1}^{n}   \big\langle  Y_i \otimes K(X_i,\cdot) , \eta \big\rangle_{\textrm{HS}}^2\\
& =  \frac{1}{n}\sum_{i=1}^{n}   \big\langle  Y_i \otimes K(X_i,\cdot) , \overline{Y \otimes K}^n - \eta \big\rangle_{\textrm{HS}}^2 \\
&\qquad+ \frac{2}{n}\sum_{i=1}^{n} \big\langle  Y_i \otimes K(X_i,\cdot) , \eta \big\rangle_{\textrm{HS}}  \big\langle  Y_i \otimes K(X_i,\cdot) , \overline{Y \otimes K}^n - \eta \big\rangle_{\textrm{HS}}
\end{aligned}
\end{equation}
we have to treat each term is this sum. First, using the  Cauchy-Schwarz inequality, the reproducing property of $K$ and condition $(\mathscr{C}_1)$, we get 
\begin{eqnarray*}
   \frac{1}{n}\sum_{i=1}^{n}   \big\langle  Y_i \otimes K(X_i,\cdot) , \overline{Y \otimes K}^n - \eta \big\rangle_{\textrm{HS}}^2 
 & \leq &\bigg( \frac{1}{n}\sum_{i=1}^{n}   \bigg\Vert  Y_i \otimes K(X_i,\cdot)\bigg\Vert ^2 _{\textrm{HS}}\, \bigg)\bigg\Vert \overline{Y \otimes K}^n - \eta  \bigg\Vert ^2 _{\textrm{HS}}\\
 & = &\bigg( \frac{1}{n}\sum_{i=1}^{n}   \left\Vert  Y_i \right\Vert ^2 _{\mathcal{Y}}\left\Vert  K(X_i,\cdot) \right\Vert ^2 _{\mathcal{H}}\, \bigg)\bigg\Vert \overline{Y \otimes K}^n - \eta  \bigg\Vert ^2 _{\textrm{HS}}\\
 & = &\bigg( \frac{1}{n}\sum_{i=1}^{n}   \left\Vert  Y_i \right\Vert ^2 _{\mathcal{Y}}   K(X_i,X_i) \, \bigg)\bigg\Vert \overline{Y \otimes K}^n - \eta  \bigg\Vert ^2 _{\textrm{HS}}\\
& \leq &\bigg( \frac{1}{n}\sum_{i=1}^{n}   \left\Vert  Y_i \right\Vert ^2 _{\mathcal{Y}}  \, \bigg)\left\Vert K\right\Vert_\infty\bigg\Vert \overline{Y \otimes K}^n - \eta  \bigg\Vert ^2 _{\textrm{HS}}.
\end{eqnarray*}
Since, from the law of large numbers, $\frac{1}{n}\sum_{i=1}^{n}   \left\Vert  Y_i \right\Vert ^2 _{\mathcal{Y}}$ and $\overline{Y \otimes K}^n $ converge in probability, as $n\rightarrow +\infty$, to  $\mathbb{E}\left(\Vert Y\Vert^2_\mathcal{Y}\right)$ and $\eta$ respectively, we deduce from the preceding inequality that 
\[
\frac{1}{n}\sum_{i=1}^{n}   \big\langle  Y_i \otimes K(X_i,\cdot) , \overline{Y \otimes K}^n - \eta \big\rangle_{\textrm{HS}}^2 =o_p(1).
\]
Secondly, using again  the  Cauchy-Schwarz inequality, the reproducing property of $K$ and condition $(\mathscr{C}_1)$, we obtain the inequality
\begin{eqnarray*}
    & & \Big| \frac{1}{n}\sum_{i=1}^{n} \big\langle  Y_i \otimes K(X_i,\cdot) , \eta \big\rangle_{\textrm{HS}}  \big\langle  Y_i \otimes K(X_i,\cdot) , \overline{Y \otimes K}^n - \eta \big\rangle_{\textrm{HS}} \Big | \\
&\leq &\bigg( \frac{1}{n}\sum_{i=1}^{n}   \bigg\Vert  Y_i \otimes K(X_i,\cdot)\bigg\Vert ^2 _{\textrm{HS}}\, \bigg)\bigg\Vert \overline{Y \otimes K}^n - \eta  \bigg\Vert  _{\textrm{HS}}\left\Vert  \eta  \right\Vert  _{\textrm{HS}}\\
& \leq &\bigg( \frac{1}{n}\sum_{i=1}^{n}   \left\Vert  Y_i \right\Vert ^2 _{\mathcal{Y}}  \, \bigg)\left\Vert K\right\Vert_\infty\bigg\Vert \overline{Y \otimes K}^n - \eta  \bigg\Vert _{\textrm{HS}}\left\Vert  \eta  \right\Vert  _{\textrm{HS}}
\end{eqnarray*}
from which we conclude that 
$$ \frac{1}{n}\sum_{i=1}^{n} \big\langle  Y_i \otimes K(X_i,\cdot) , \eta \big\rangle_{\textrm{HS}}  \big\langle  Y_i \otimes K(X_i,\cdot) , \overline{Y \otimes K}^n - \eta \big\rangle_{\textrm{HS}}= o_p(1). $$
Consequently, from \eqref{diff} it is seen that $ \frac{1}{n}\sum_{i=1}^{n}   \big\langle  Y_i \otimes K(X_i,\cdot) , \overline{Y \otimes K}^n \big\rangle_{\textrm{HS}}^2$ has the same limit in probability than $\frac{1}{n}\sum_{i=1}^{n}   \big\langle  Y_i \otimes K(X_i,\cdot) , \eta \big\rangle_{\textrm{HS}}^2$. From the law of large numbers this latter converges in probability, as $n\rightarrow +\infty$ to     $\mathbb{E} \left(\big\langle  Y_1 \otimes K(X_1,\cdot) , \eta \big\rangle_{\textrm{HS}}^2\right)$.
On the other hand, we have the inequality
\begin{eqnarray*}
& &\Big\vert \frac{1}{n}\sum_{i=1}^{n}   \big\langle  Y_i \otimes K(X_i,\cdot) , \overline{Y \otimes K}^n \big\rangle_{\textrm{HS}}    - \frac{1}{n}\sum_{i=1}^{n}   \big\langle  Y_i \otimes K(X_i,\cdot) , \eta \big\rangle_{\textrm{HS}} \Big\vert  \\
& =&\Big\vert   \frac{1}{n}\sum_{i=1}^{n}   \big\langle  Y_i \otimes K(X_i,\cdot) , \overline{Y \otimes K}^n - \eta \big\rangle_{\textrm{HS}}\Big\vert \\
& \leq &\bigg( \frac{1}{n}\sum_{i=1}^{n}   \bigg\Vert  Y_i \otimes K(X_i,\cdot)\bigg\Vert _{\textrm{HS}}\, \bigg)\bigg\Vert \overline{Y \otimes K}^n - \eta  \bigg\Vert _{\textrm{HS}}\\
 & = &\bigg( \frac{1}{n}\sum_{i=1}^{n}   \left\Vert  Y_i \right\Vert  _{\mathcal{Y}}\left\Vert  K(X_i,\cdot) \right\Vert _{\mathcal{H}}\, \bigg)\bigg\Vert \overline{Y \otimes K}^n - \eta  \bigg\Vert  _{\textrm{HS}}\\
 & = &\bigg( \frac{1}{n}\sum_{i=1}^{n}   \left\Vert  Y_i \right\Vert _{\mathcal{Y}}   \sqrt{K(X_i,X_i)} \, \bigg)\bigg\Vert \overline{Y \otimes K}^n - \eta  \bigg\Vert _{\textrm{HS}}\\
& \leq &\bigg( \frac{1}{n}\sum_{i=1}^{n}   \left\Vert  Y_i \right\Vert  _{\mathcal{Y}}  \, \bigg)\left\Vert K\right\Vert_\infty^{1/2}\bigg\Vert \overline{Y \otimes K}^n - \eta  \bigg\Vert  _{\textrm{HS}}
\end{eqnarray*}
which implies that  $n^{-1}\sum_{i=1}^{n}   \big\langle  Y_i \otimes K(X_i,\cdot) , \overline{Y \otimes K}^n \big\rangle_{\textrm{HS}}    - n^{-1}\sum_{i=1}^{n}   \big\langle  Y_i \otimes K(X_i,\cdot) , \eta \big\rangle_{\textrm{HS}} =o_p(1)$ since, from the law of large numbers,  $\frac{1}{n}\sum_{i=1}^{n}   \left\Vert  Y_i \right\Vert  _{\mathcal{Y}}$ and $\overline{Y \otimes K}^n $ converge in probability, as $n\rightarrow +\infty$, to  $\mathbb{E}\left(\Vert Y\Vert_\mathcal{Y}\right)$ and $\eta$  respectively. Consequently,  $n^{-1}\sum_{i=1}^{n}   \big\langle  Y_i \otimes K(X_i,\cdot) , \overline{Y \otimes K}^n \big\rangle_{\textrm{HS}}$ converges in probability, as $n\rightarrow +\infty$,  to the same limlit than  $n^{-1}\sum_{i=1}^{n}   \big\langle  Y_i \otimes K(X_i,\cdot) , \eta \big\rangle_{\textrm{HS}}$ , that is $\mathbb{E}\left(\big\langle Y_1 \otimes K(X_1,\cdot) , \eta \big\rangle_{\textrm{HS}}\right)$. Finally, from \eqref{alpha}, we deduce that   $\widehat{\alpha}^2$ converges in probability, as $n\rightarrow +\infty$, to  $Var\left( \big\langle  Y_1 \otimes K(X_1,\cdot),  \eta \big\rangle_{\textrm{HS}}\right)$.  
%\section*{References}

\end{document}